# Bio-Heat Transfer and Monte Carlo Measurement of Near-Infrared Transcranial Stimulation of Human Brain

**Faezeh Ibrahimi[1] and Mehdi Delrobaei[1,*]**

[1]The Center for Research and Technology (CREATECH), Faculty of Electrical Engineering, K. N. Toosi University of Technology, Tehran, Iran
*delrobaei@kntu.ac.ir

## ABSTRACT

Transcranial photobiomodulation is an optical method for non-invasive brain stimulation. The method projects red and near-infrared light through the scalp within 600-1100 nm and low energy within the 1-20 $J/cm^2$ range. Recent studies have been optimistic about replacing this method with pharmacotherapy and invasive brain stimulation. However, concerns and ambiguities exist regarding the light penetration depth and possible thermal side effects. While the literature survey indicates that the skin temperature rises after experimental optical brain stimulation, inadequate evidence supports a safe increase in temperature or the amount of light penetration in the cortex. Therefore, we aimed to conduct a comprehensive study on the heat transfer of near-infrared stimulation for the human brain. Our research considers the transcranial photobiomodulation over the human brain model by projecting 810 nm light with 100 $mW/cm^2$ power density to evaluate its thermal and optical effects using bioheat transfer and radiative transfer equation. Our results confirm that the near-infrared light spectrum has a small incremental impact on temperature and approximately penetrates 1 cm, reaching the cortex. A time-variant study of the heat transfer was also carried out to measure the temperature changes during optical stimulation.

## Introduction

In the last two decades, efforts to find an efficient yet non-invasive cognitive rehabilitation method resulted in the development of the therapeutic application of optical brain stimulation, known as photobiomodulation (PBM) therapy[1,2]. Photobiomodulation applies red and near-infrared light within 600-1100 nm and low energy within the 1-20 $J/cm^2$ range to penetrate and stimulate the tissues[3]. Some beneficial effects attributed to this method such as: (1) the CcO enzyme (Cytochrome c Oxidase) stimulation in the mitochondrial respiratory chain[4], (2) activation of transcription factors, (3) induction of gene expression by production and distribution of $Ca^{2+}$ ions, (4) improving neural metabolic capacity, (5) preventing neuronal apoptosis, and (6) and repair and re-growth of neural synapses[3,5-9].

Various clinical studies have been conducted to evaluate the light-based brain stimulation methods. The most effective wavelengths in reducing cognitive disease symptoms were 810 nm, 850 nm, 670 nm, and 710 nm[3]. The 810 nm light has shown the deepest penetration. This particular light is currently being used in commercial brain stimulation devices[10]. However, implemented investigations are more qualitative, lacking detailed analysis of 810 nm light penetration depth in the brain tissue, particularly in the cortex layer. Furthermore, the photothermal effect has also raised concerns[11,12].

The PBM therapy still faces two critical challenges: (1) light diffusion into the brain tissue[3] and (2) photothermal side effects[11]. A literature survey indicates no consistent evidence for near-infrared depth penetration. The experimental and simulation results on the side effects of PBM seem to be inconsistent. About 1% of the fluence rate of 670 nm light reaches the cerebral cortex of the rat[13]. This number is estimated to increase to 2% for 810 nm in the human brain[12]. These results have not been fully validated, and the research is still ongoing.

Therefore, we aimed to conduct a comprehensive simulation study using a manually assembled high-resolution geometry and a fine mesh with more light sources on the scalp. We considered an air gap between the light sources and the skin to more accurately evaluate the penetration depth and the thermal effects of the 810 nm light. Thermal and optical efficacy on the neural tissues were estimated using bioheat transfer and radiative transfer equations.

## Methods

We initially assembled a complete geometry, refined the mesh, defined each material properties with high exactitude based on different sources, and tried to comply with clinical investigations. Hence, our simulation method consists of geometry and material selection, physics, boundary conditions, and meshing. Eventually, the results are fine-tuned with an impedance of node and validated using the results already presented in the literature. The following subsections describe the simulation steps and define the employed equations and the boundary conditions in COMSOL Multiphysics software, version 6.0 (COMSOL Inc. Stockholm, Sweden).

**Geometry Selection**

The proposed human head model presented in this work consists of three separate parts: the brain, skull, and scalp. The scalp model is provided by a standard model in the COMSOL Multiphysics source and used as the standard base for the other components[14]. The skull model is a 3D mesh based on MRI images and modified by a virtual geometry tool to remove the sharp edges and merge the faces. The brain's model is a standard model based on high-resolution MRI images imported to the COMSOL environment. We planned to cover different parts of the cerebral cortex, similar to the experimental brain stimulation reports using a cap with a set of light-emitting diodes (LEDs). Hence, we considered nine LEDs equally distributed on the head surface. The LEDs were modeled based on the standard 810 nm LEDs already used in PBM experiments. The advantages of using such LEDs include: (1) appropriate surface-to-volume ratio (compactness), (2) easy installation, (3) effective thermal insulation, and (4) low energy consumption (high efficiency). The final model forms twelve separate components inside a cylindrical air configuration. Fig. 1 represents the final geometry components for this study.

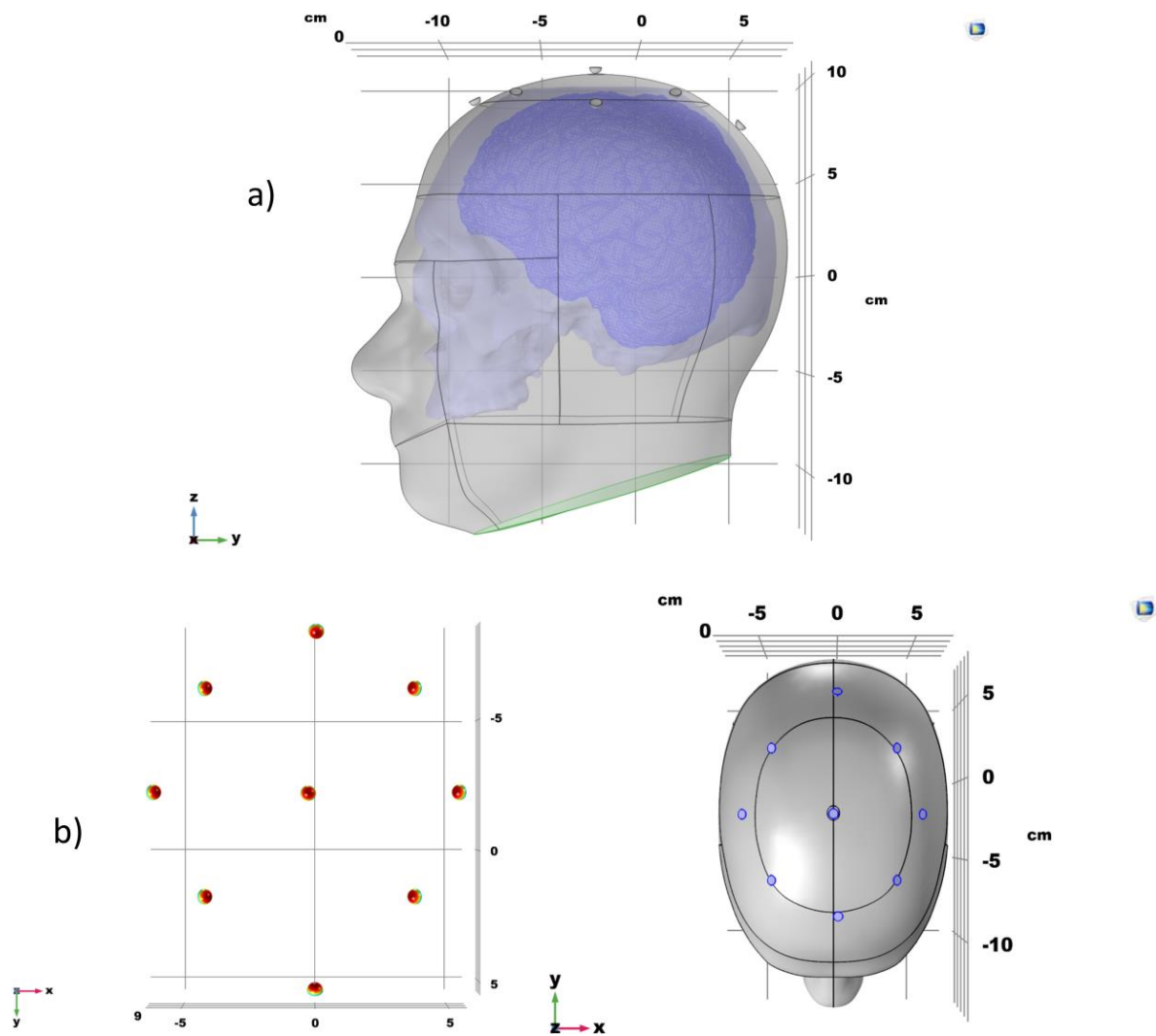

**Figure 1.** (a) The final assembled model, (b) Location of nine LEDs on the skin.

## Parameter Selection

In this step, we investigate the properties of the tissues and the LEDs. Table 1 lists the parameters provided for the heat transfer equation. In this table, $K$, $C_p$, and $\rho$ are the conduction heat transfer, the specific heat capacity in the constant pressure, and the mass density coefficients, respectively.

**Table 1.** Properties of the material for the whole components[15–17]

| Properties Components | K $(W/m.°C)$ | $\rho$ $(kg/m^3)$ | $C_p$ $(J/kg.°C)$ |
|---|---|---|---|
| Scalp | 0.50 | 1200 | 4000 |
| Skull | 1.15 | 1990 | 2300 |
| Brain | 0.57 | 1050 | 3650 |
| LED | 0.20 | 1190 | 1170 |

The mentioned information is inadequate for the whole bioimaging study as we must consider the bioheat transfer mechanism and optical photon transport concepts. Therefore, more variables will be defined as partial differential equations.

## Physics and the Boundary Conditions

Multiphysics implemented in this simulation includes: (1) a coefficient form of a partial differential equation (PDE) to model the radiative transfer equation to the fluence rate calculation, and (2) the bioheat transfer physic and measures the temperature variation in the tissues under the heat light influence.

### Diffusion Approximation through PDE form Physics in RTE Simulation

A diffusion approximation of the radiative transfer equation (RTE) computes the fluence rate through the following second-order PDE[18,19,22].

$$\left(\frac{1}{v}\frac{\partial}{\partial t} + \hat{s}\nabla + \mu_a(r) + u_s'(r)\right) L(r,\hat{s},t) = Q(r,\hat{s},t) + \mu_t(r) \int f(s,\hat{s},r) L(r,\hat{s},t)^2 \hat{s} \tag{1}$$

where $\mu_a$, $u_s'$, and $\mu_t$ are the absorption, reduced scattering, and total attenuation coefficients, and v refers to the velocity of the light through the tissues. The term $L(r,\hat{s},t)$ is the radiance at position r with the direction of propagation S and $Q(r,\hat{s},t)$ defined as the source term.

The following expressions describe isotropic fluence rate, $\phi(r,t)$, and a small directional flux, $J(r,t)$:

$$\phi(r,t) = \iint_{4\pi} L(r,\hat{s},t) w \tag{2}$$

$$J(r,t) = \iint_{4\pi} L(r,\hat{s},t) \hat{s} w \tag{3}$$

Eventually, the final diffusion equation can be described as:

$$\frac{1}{v}\frac{\partial \phi}{\partial t} - \nabla D(r) \nabla \phi(r,t) + \mu_a(r) \phi(r,t) = Q_0(r,t) \tag{4}$$

where D(r) and $\mu_s'$ are derived as:

$$D(r) = \frac{1}{3(\mu_a(r) + \mu_s'(r))} \tag{5}$$

$$\mu_s'(r) = \mu_s(1-g) \quad ; \quad g = 0.89 \tag{6}$$

$$V = c/n \tag{7}$$

**Table 2.** The optical properties of the biological tissues in a human head model at 810 nm[23]

| Properties Tissues | $\mu_a$ (1/cm) | $\mu_s$ (1/cm) | D (cm) |
|---|---|---|---|
| Scalp and skull | 0.16 | 0.76 | 0.043 |
| Brain | 0.57 | 0.80 | 0.039 |

Equation 4 utilizes the PDE form physics with calculated D and $\mu'_s$ from Table 2 for each layer (components and chromophores were efficient in calculating $\mu_a$ of the brain.). Consider c, the velocity of light in vacuum, and n, the refractive index.

For evaluating the fluence rate $\varphi(r, t)$, two separate equations following Eq. 4 were employed for the tissues based on Table 2. We assumed n = 1.4 and the anisotropy factor, g = 0.89 for all vital layers. It is noted that part of the energy consumed by a light source is converted to light and the rest of the energy to heat. For a standard LED and a fluence rate of 100 $mW/cm^2$, we assumed that 65% of the energy is converted to heat and 35% to light. Therefore, 65 $mW/cm^2$ has been utilized in Eq. 6 and 35 $mW/cm^2$ in Eq. 4.

The Dirichlet boundary condition was applied for the backside of LEDs and the human scalp's surfaces, which are not on the light path. The sources were placed at the air-tissue interface close to the scalp following Chaieb et al.[24].

### *Bioheat Transfer Physics for Evaluating Thermal Effect of Photobiomodulation*

The following bioheat transfer equation was utilized based on Penne's bioheat equation[25,26].

$$\rho C_p \frac{d}{dt} \nabla T - K \nabla^2 T = q_s + [\rho_b c_{p,b} w_b (T_b - T) + q_{met}] \qquad (8)$$

where $\rho$, $C_p$, and K are density, specific heat capacity, and conduction heat transfer coefficients for tissues, and $q_s$ depends on the heat source. The rest of the coefficients are related to the blood flow; for instance, $w_b$, $T_b$, and $q_{met}$ are blood perfusion, temperature, and metabolic heat coefficients.

This relationship integrates the heat transfer mechanism with the fluid effect of cerebral blood flow. Our light sources generate heat by assuming 65% heat loss in the LEDs. Therefore, 65 $mW/cm^2$ was assigned in a boundary heat source for each LED. The initial temperature for the brain and the skull was set to 37 °C. The initial temperature for the skin was set to 33 °C. The LEDs were given an initial temperature of 25 °C, the same temperature as the environment.

As a boundary condition, it is assumed that the heat loss from the scalp (to the environment) occurs due to convection (the air surrounding the model is considered stationary with the convection coefficient of h= 0.5 $mW/cm^2$.°C and the temperature Tinf = 25 °$C^{27}$). The same boundary condition was used for the LED's external surface.

The bioheat transfer includes the blood flow, density, temperature, volumetric metabolism, and specific heat for the human blood in a neutral condition. Table 3 incorporates this information considering the blood density and the blood specific heat to be 1050 $kg/m^3$ and 3600 $J/kg.°C$, respectively. In the following steps, we divide the geometry for the finite element method

**Table 3.** Bioheat transfer properties considering blood flow[12]

| Properties Tissues | Metabolic Heat (W/m³) | Blood Perfusion (1/s) |
|---|---|---|
| Scalp | 363 | 0.00143 |
| Skull | 70 | 0.000143 |
| Brain | 10437 | 0.08 |

and choose how to study and verify the results.

### **Meshing**

It is noted that the meshing process has a significant role in the speed and accuracy of the computations. Therefore, the small partitions' size, type, and distribution are essential for accurate results. The model's meshing process details can be seen in Table 4. The brain's model requires small and compressed elements to more accurately reflect the principal influence of the PBM on the cortex layer. The degrees of freedom were solved for 12,402,370 domain elements, 1,246,350 boundary elements, and 653,949 edge elements. The whole geometry has been provided as supplemental material (available online).

Table 4. Mesh parameters at optimized physics control mesh for each component

| Components | Size | Type | Number of Domain elements | Number of Boundary elements | Number of Edge elements |
|---|---|---|---|---|---|
| Head | Normal | Tetrahedral | 155921 | 29360 | 463 |
| Skull | Fine | Tetrahedral | 4476531 | 1231562 | 652992 |
| Brain | Finer | Tetrahedral | 6805110 | 1099348 | 629762 |
| LED | Finer | Tetrahedral | 3097 | 1570 | 216 |

## Results

The results presented in this section are based on the defined hierarchy in the stationary process. Fig. 2 depicts a suggested cutline in the z-axis to analyze variables in the form of a line graph. Cutline is a parameter analysis path to evaluate variation in our software.

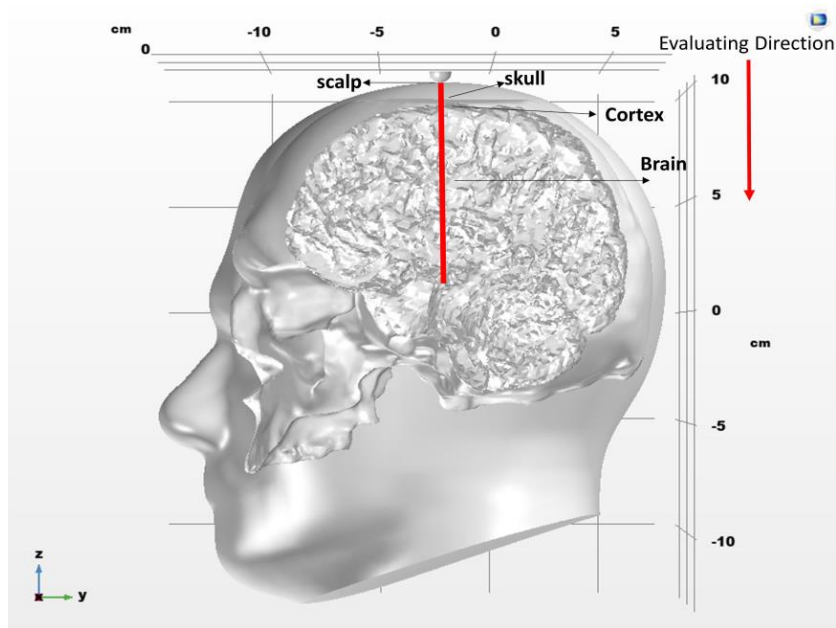

**Figure 2.** The cutline parallel to the z-axis.

### Thermal effect of Phtotobiomodulation at 810 nm

Fig. 3(a). represents a 3D contour of the steady-state temperature for the whole tissues affected by PBM with 810 nm light and 65% heat loss. A utilized color spectrum describes a higher temperature in $°C$ with a dimmer color. Although each tissue color seems uniform, slight differences exist in several sections of every layer. The temperature contour obtained from bioheat transfer physics, boundary conditions, and bioheat transfer properties.

### Temperature distribution on the cutline path

Fig. 3(b). depicts temperature variation in Celsius ($°C$). Note that the origin of coordinate is located in the brain with almost 10 cm depth under the skin and the whole temperature differences are less than 0.5 degrees. Line graph showing temperature from 37.02 $°C$ to 37.45 $°C$ on the y-axis against z-coordinate within 2 cm on the x-axis. According to the following figure, the highest temperature in tissues during PBM occurs in the scalp at almost 37.45 $°C$. Temperature declines between the skull and brain in the cortex layer and will set to approximately 37.035 $°C$ for the brain. Due to the line graph, the scalp temperature is varied from 37.15 to 37.45 $°C$ and has a higher temperature than the blood temperature below the skin. The dash lines draw roughly to separate the tissue boundaries.

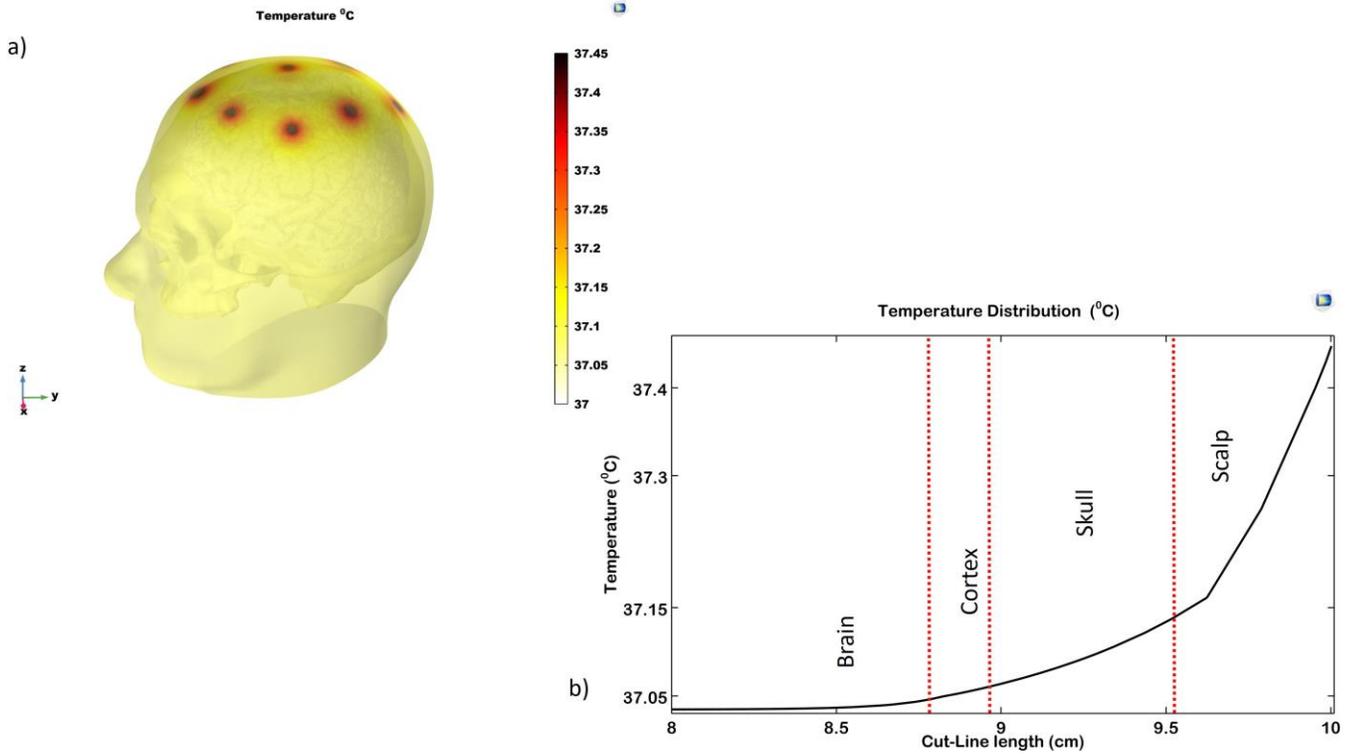

**Figure 3.** (a) Temperature contour on the scalp of the head model., (b) Temperature variation in whole tissues affected by 810 nm light.

### Heat absorbed by the brain through photobiomodulation on the cutline path

The absorbed heat flux must be studied for thermal damage in every component. Absorbed heat flux is calculated based on the tissue's conduction heat transfer coefficient and temperature variation. The volume of heat flux generated in the brain tissue is calculated as follows:

$$Q = mc_p \Delta T \quad \text{and} \quad m = \rho v \quad \text{then} \quad heat\ flux = \frac{Q}{v} = \rho c_p (T - T_i) \quad (9)$$

where m, $c_p$, and $\Delta T$ are the mass, specific heat capacity in the brain tissue's constant pressure and temperature variation. Terms $\rho$ and v refer to the density and volume, and $T_i$ is defined as the brain's initial temperature. It is worth mentioning that the term T is the final temperature, and the heat flux unit is $\frac{J}{cm^3}$. (We considered $\rho = 0.00105\ \frac{kg}{cm^3}$, $C_p = 3650\ \frac{J}{kg.°C}$, and $T_i = 37\ °C$.)

Fig. 4(a). indicates the brain's heat absorbed per unit area. The absorbed heat flux has been assumed to be generated by the temperature alteration before and during the PBM. According to the results, volume heat flux is variable in the cortex domain within 0.12-0.21 $J/cm^3$. However, the rest of the brain tissue has a constant value of fewer than 0.15 $J/cm^3$. Note that the cortex layer thickness is considered 2.5 mm on average.

### Fluence rate for the whole tissues through Phtotobiomodulation at 810 nm

The following line graph in Fig. 4(b). represents the fluence rate value by $mW/cm^2$ in this case study. The fluence rate for layers is 0-0.36 $mW/cm^2$. Our results are based on an approximately RTE simulation equation, and absorption and scattering are considered for each tissue. According to Fig. 4(b)., the light energy is dissipated, passing through the upper tissues. Due to the light attenuation, the fluence rate in brain tissue is fewer than in the skull. Light attenuation occurs in the cortex layer of the brain, and part of it is affected by the slight amount of light. Note that the 810 nm light can penetrate up to a depth of 1 cm.

### Temperature and Fluence rate contour for the brain tissue

In the outcomes of the preceding figures, insignificant differences exist in each tissue. General contours were exhibited to compare tissues with each other. In this section, temperature and fluence rate contour are plotted only for the brain tissue. The ultimate aim of this imaging is to examine optically and thermally measure the surface of the cerebral cortex. Fig. 5 indicates temperature variation in section a, and the fluence rate contours in section b for the brain tissue. These results confirm insignificant temperature rise (max = 0.05 °C variation) and cortex optical stimulation through 810 nm transcranial PBM. According to Fig. 5 section b, nine LEDs layouts covered the cerebral cortex, and all sources could bring light to the brain's surface, even in small amounts.

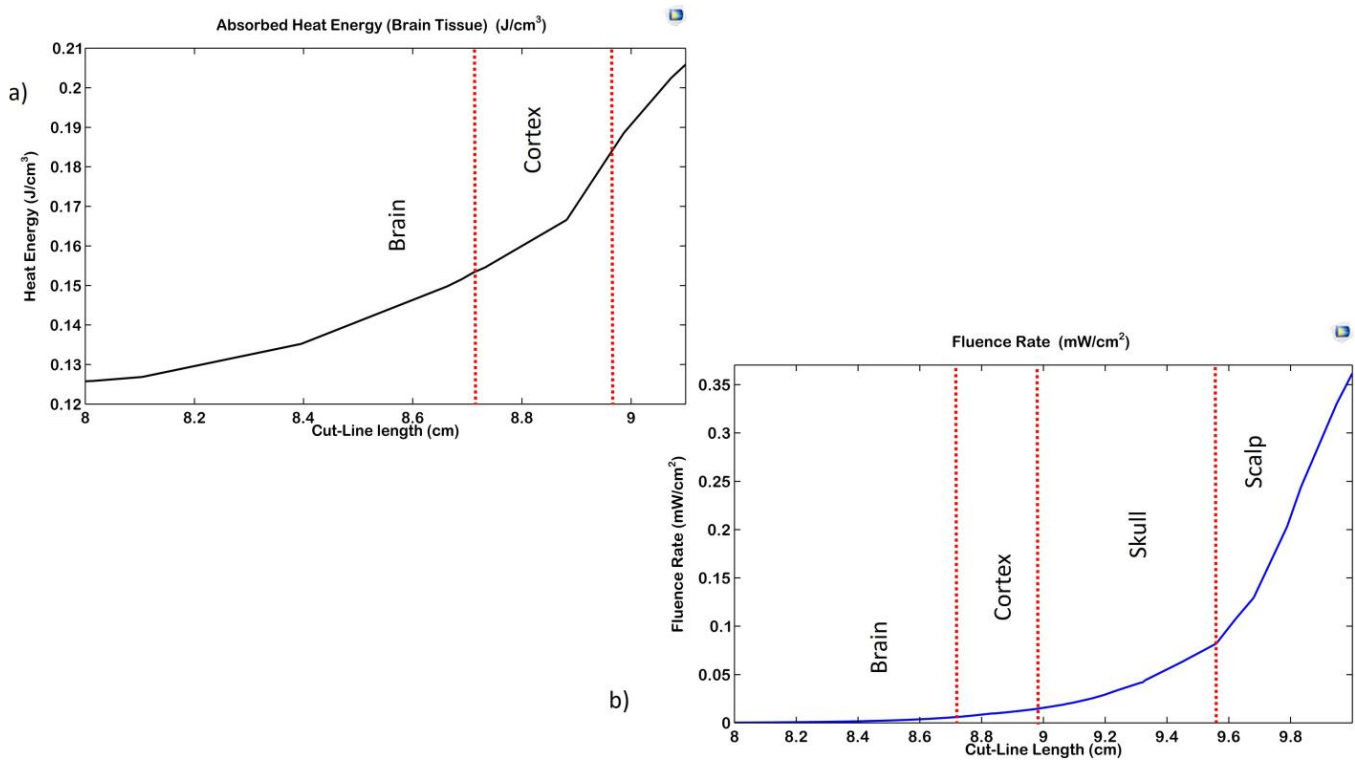

**Figure 4.** (a) Heat absorbed by the brain through 810 nm PBM, (b) Fluence rate variation in whole tissues on the cutline path.

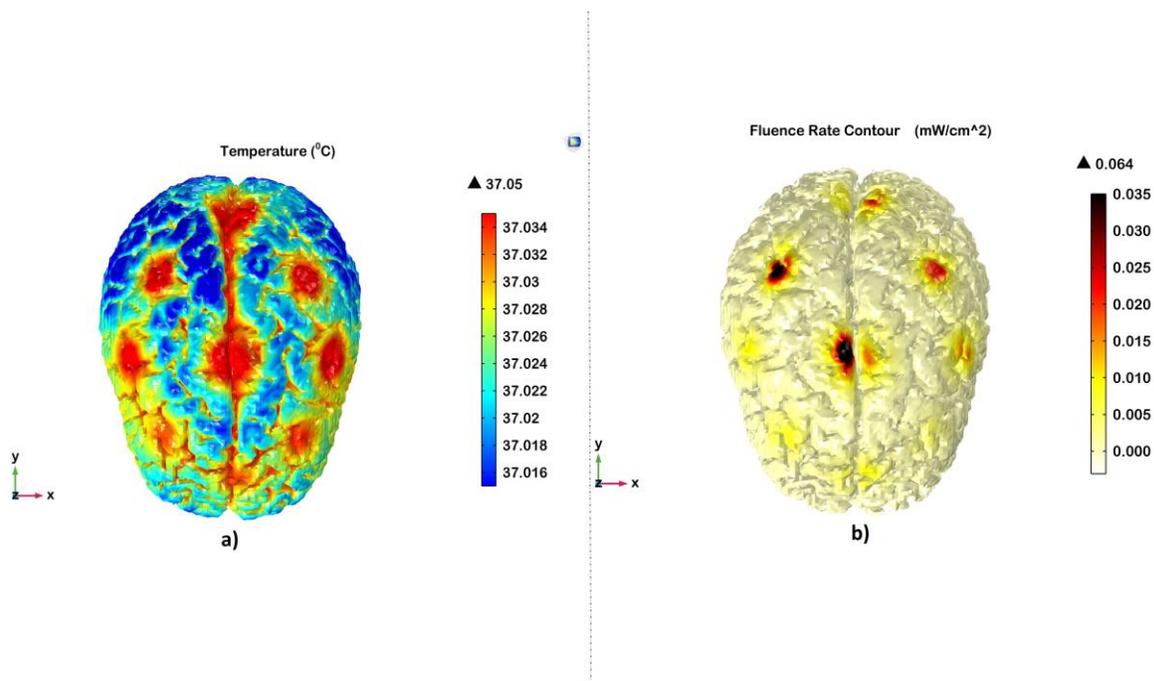

**Figure 5.** (a) Temperature contour; (b) Fluence rate contour for the brain tissue.

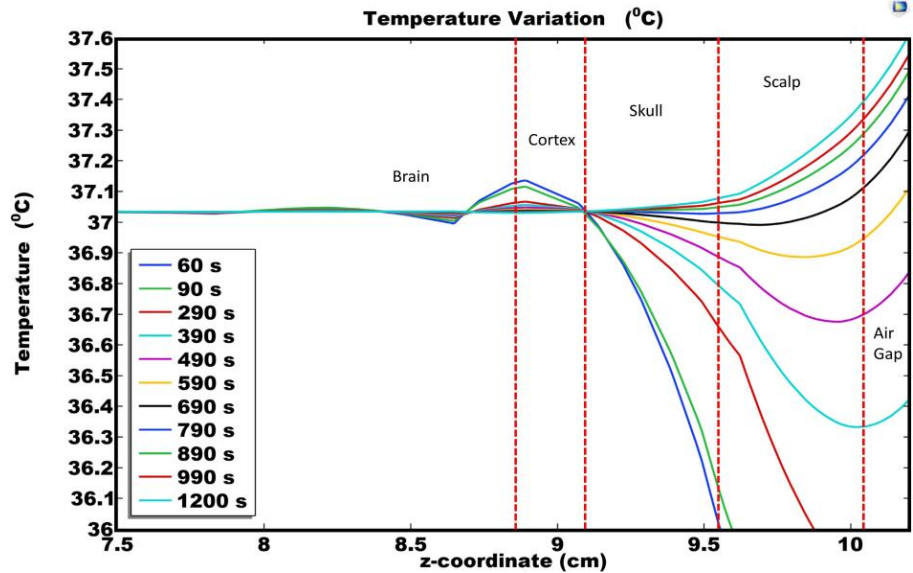

**Figure 6.** Steady-state analysis graph for temperature through all tissues during NIRS stimulation.

### Temperature Time-Variant Analysis

A time-variant study was conducted with the same physics and conditions within 20 minutes (1200 s), and the results are represented in Fig. 6. This analysis helps to understand the temperature variation during brain stimulation. According to Fig. 6, the temperature fluctuates due to the 810 nm transcranial stimulation, considering the scalp has a maximum change (almost 1.6 °C). It is noted that changes in the brain tissue's temperature are minor. Over time, the temperature curve for all the tissues becomes closer to the steady-state result in Fig. 3. Furthermore, after the stimulation period of about 1200 s, they conform without thermal side-effect.

### Sensitivity Analysis

A sensitivity analysis of the input parameters was performed to investigate the fluence rate variations. The input parameters include wavelength and the power of the input light sources. Fig. 7(a) represents the rate of change in the fluence rate based on the energy of the light source. As predicted, higher light energies have a higher fluence rate, hence more penetration. In Fig. 7(b), three operational wavelengths in brain stimulation are compared, and we observed that the most effective wavelength is 810 nm, consistent with the results of clinical trials.

## Discussion

Researchers have raised concerns that the penetration depth and the thermal effects of optical brain stimulation (as a side-effect) are not yet fully known[3,11]. The proposed work tries to address this unmet need.

Considering thermal results for skin temperature in the prior simulation study[12] and comparison with the experimental results[11] exhorted us to investigate further. Therefore, we endeavored to repeat this survey in a different geometry with a new refined mesh and input power conversion ratio into light and heat for each source (65% heat loss and 35% light). Results can be different depending on geometry, physics, and boundary conditions.

We collected the geometry components based on the standard human head model and MRI images. Meshing was implemented with high accuracy, mainly in the brain. The principle of node independence was considered. It is noted that we used high-performance computing for the simulation study.

The clinical studies and simulations are complementary and should be used in a different context for validation. For example, skin temperature is known in experimental research and is used as the basis for simulation results. In comparison, internal tissues' temperature estimating, such as the brain, is obtained from computational results. Our results show conformity in both aspects.

It is noted that most PBM studies are clinical, and there is little research with simulation content in this field. The beneficial effects of using PBM as a therapy were studied in most experimental research, and its challenges had been less discussed.

To acquire more accurate results and higher computational efficiency in the fluence rate for complex media, such as biological tissues investigation, the Radiative Transport Equation (RTE) essentially is used[19,20,21,22]. The proposed model was a low-absorption high-scattering medium, and a diffusion approximation of the RTE (Eq. (4)) was used for comparison, mainly due to its rapid and relatively reliable solution. Note that diffusion approximation is a method for solving RTE, and considering the fast response, the light energy will be damped faster than the other methods. Therefore, the reported data proposes the minimum amount of fluence rate for specific cortex tissue.

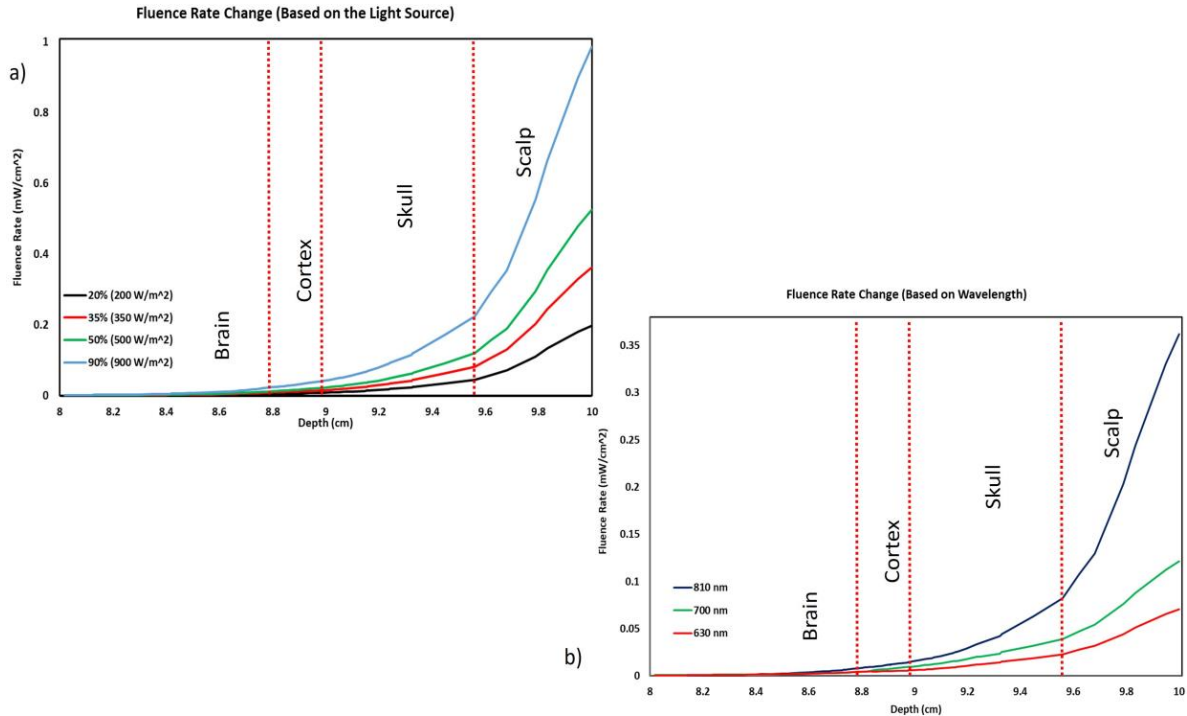

**Figure 7.** Sensitivity analysis graphs for all tissues based on (a) input light source and (b) wavelength.

Temperature variation induced by the light absorption was estimated using bioheat transfer physics based on Penn's equation[26]. As predicted, the maximum value of temperature refers to the scalp. Because the skin has a heat loss through convection heat transfer with the environment but is located near heat sources and the brain has a more distance from the 100 $mW/cm^2$ sources[28] and the blood perfusion act as a heat sink for the brain tissue, which considered.

A more detailed investigation illustrated in Fig. 3 depicted that the skin temperature variances within 37.15- 37.45 °C affected by 810 nm PBM will not cause burning or any unpleasant impression. This temperature scale is consistent with previous simulations and clinical reports. According to Fig. 3, for the cortex layer with a 2.5 mm average thickness, the calculated temperature rise was less than 0.05 °C. This value is also about 0.045 °C in a similar work[12].

Fig. 5 confirmed approximately 0.035 °C increase in the cortex layer and proved the 810 nm light reaches this layer. Also, nine LEDs arrangement on the cerebral cortex has acceptably covered the brain's surface. However, this number can vary. We only tried to provide a better representation of optical brain stimulation, as shown in Fig.5, given that clinical studies are often performed with more LEDs.

Fig. 4 represents the volume of heat absorbed by the brain tissue, for the cortex layer was a negligible amount within 0.12-0.2 $J/cm^3$, which is safe concerning thermal effects. The fluence rate variation is represented in Fig. 4 to evaluate the sufficiency of penetrated 810 nm light into tissue depth. The proposed results prove partial penetration into the cortex layer for efficient depth (effective light penetration is about 1 cm). However, according to the experimental research, this amount of 810 nm light showed a significant therapeutic effect on cognitive and neurodegenerative diseases[3,28].

As changes in temperature, even less than 1 °C, affect the brain's cognitive function due to the high thermal sensitivity of its neurons, temperature control is a critical challenge for transcranial PBM. According to Fig. 6, during 20 minutes of stimulation, the brain temperature does not exceed the allowable limit. Moreover, we confirmed that reported results (Fig. 3- 5) are available after at least 20 minutes of stimulation. It is noted that the time for clinical brain stimulation is about 20- 30 minutes.

Within the first minute, we observed a rapid temperature increase, following the changes from initial values for every tissue per distance from the sources and physical properties. This sharp turn and temperature increase for the cortex was less than 0.15 °C, which seems to be a logical value. Over time, the effect of blood flow was more substantial than the received heat energy and modulated this temperature increase until reaching a steady state. According to Fig. 6, none of the simulated temperatures were in a dangerous range.

A comprehensive simulation study reported a 0.11 $J/cm^2$ fluence value in the human brain for 300 $mW/cm^2$ in 5 minutes of 810 nm PBM therapy[29]. This result has good conformity with Fig. 4(b) and Fig. 5(b), as we reported 0.036 $mW/cm^2$ fluence rate for 100 $mW/cm^2$ light source in 20 minutes (almost 0.045 $J/cm^2$ fluence). Another research reported a 0.12 °C temperature increase in cortex tissue for 808 nm laser with an intensity of 318 $mW/cm^2$ in 5 minutes[30]. This result complies highly with our result shown in Fig. 6 where we reach this value after 2 minutes. However, the steady-state temperature for the cortex (after 20 minutes) is less than 0.05 °C

Our efforts in this study focused on developing a comprehensive human head model and proposing accurate simulation results; however, some limitations exist. For instance, we did not compare the wavelengths that reduce neural damage and measure the most effective wavelength. We also did not model the cerebral blood flow (CSF) layer. Considering more brain layers may increase the accuracy of the results and improve compliance with clinical outcomes.

This study aimed to clarify two critical issues in optical brain stimulation through 810 nm PBM. First, to determine if adequate light energy with a specific wavelength reaches the brain cortex. Second, to investigate the level of temperature rise, which can be harmful.

Gagnon et al. listed the optical properties of skin, skull, cerebrospinal fluid (CSF), and white and gray matter of the brain tissue[31]. It is evident that the CSF's absorption and scattering coefficients (1.7 1/$m$ and 10 1/$m$, respectively) are pretty small compared to other layered tissues (for instance, 16.4 1/$m$ and 740 1/$m$ for the skin). Shimojo et al. also reported the absorption and reduced scattering coefficients for the blood[32], which is more effective compared to the CSF layer for 810 nm wavelength[12]. Another research concluded that the existence of the CSF layer does influence the photon pathlengths in other tissue layers but contributes little to the absorption and scattering of light and can be ignored[33].

In our future work, we plan to expand the analysis by considering more layers of the brain (for instance, CSF), extending the analysis time to more than 30 minutes, and having a spherical investigation for the temperature and fluence rate over time.

## Data availability

All data are available in the main text or supplementary information.

## Acknowledgments

We thank the Cognitive Science and Technologies Council (CoGC) of Iran for supporting this study.

## Author contributions statement

Conceptualization, F.I. and M.D.; Methodology, F.I. and M.D.; Investigation, F.I.; Software, F.I.; Formal Analysis, F.I.; Writing–Original Draft, F.I.; Review Editing, F.I. and M.D.; Supervision, M.D.

## Competing interests

The authors declare no competing interests.

## Additional information

Supplementary Material The complete geometry is available through the link provided below. The proposed geometry is in STL and MPHBIN formats, and the model is surrounded by air. All the tissues are available and assembled, and it is recommended to use the transparency tool to acquire a better observation. `https://drive.google.com/file/d/1nyuSErvOneidjuI80lgw8uR7yY5XzH5-/view?usp=sharing`